\documentclass[12pt]{article}
\usepackage{amssymb, amsmath}

\newtheorem{theorem}{Theorem}[section]

\newcommand{\n}{\nonumber}

\newcommand{\si}{\sigma_R }

\newcommand{\s}{\sigma}

\newcommand{\bb}{\begin{equation}}
\newcommand{\ee}{\end{equation}}
\newcommand{\bq}{\begin{eqnarray}}
\newcommand{\eq}{\end{eqnarray}}
\newcommand{\bqn}{\begin{eqnarray*}}
\newcommand{\eqn}{\end{eqnarray*}}

\begin{document}
\title{ On the Liouville type theorem for stationary compressible Navier-Stokes-Poisson equations
in $\Bbb R^N$}
\author{ Dongho Chae\\
 Department of Mathematics\\
  Sungkyunkwan
University\\
 Suwon 440-746, Korea\\
 e-mail: {\it chae@skku.edu }}
 \date{}
\maketitle
\begin{abstract}
In this paper we prove Liouville type result for the stationary
solutions to the compressible Navier-Stokes-Poisson equations(NSP)
and the compressible Navier-Stokes equations(NS)  in $\Bbb R^N$,
$N\geq 2$. Assuming suitable integrability and the uniform
boundedness conditions for the solutions we are led to the
conclusion that $v=0$. In the case of (NS) we deduce that the
similar integrability conditions imply $v=0$ and $\rho=$constant on
$\Bbb R^N$. This shows that if we impose the the non-vacuum boundary
condition at spatial infinity for (NS), $v\to 0$ and $\rho\to
\rho_\infty
>0$, then $v=0$, $\rho=\rho_\infty$ are the solutions.\\
\ \\
{\bf AMS Subject Classification Number:}76N10,76N15, 35Q35 \\
{\bf keywords:}  stationary compressible Navier-Stokes Poisson
equations, stationary compressible Navier-Stokes equations,
Liouville type theorem
\end{abstract}
\section{Introduction}
 \setcounter{equation}{0}
We are concerned  on the  stationary Navier-Stokes-Poisson equations
in $\Bbb R^N$, $N\geq 2$.
  \bq
&&\label{11} \mbox{div}(\rho v) = 0, \\
&&\label{12}   \mbox{div}(\rho v \otimes v) = -\nabla
p+k\rho \nabla \Phi+ \mu \Delta v +(\mu +\lambda ) \nabla \mathrm{div }\, v, \\
&\label{13} &\Delta \Phi =\rho-\rho_0, \\
 &&\label{14}  p=a \rho^\gamma, \, \gamma >1.
 \eq
The system (NSP) describes  compressible gas flows,  and $\rho, v,
\Phi$ and $p$ denote the density, velocity, the potential of the
underlying force and the pressure respectively.  The constant
$\rho_0 \geq 0$ is called the background state. The viscosity
constants $\lambda, \mu$ satisfy $\mu>0$, $\lambda +\mu
>0$. Here $k$ is also a physical constant, which signifies the
property of the forcing, repulsive if $k > 0$ and attractive if $k <
0$. In particular, in the case $k=0$( and without (\ref{13})) the
system reduces to the following stationary compressible
Navier-Stokes equations(NS).
 \bq
&&\label{11a} \mbox{div}(\rho v) = 0, \\
&&\label{12a}   \mbox{div}(\rho v \otimes v) = -\nabla
p+ \mu \Delta v +(\mu +\lambda ) \nabla \mathrm{div }\, v,  \\
 &&\label{14a}  p=a \rho^\gamma, \, \gamma >1.
 \eq
 For the
general treatment of the system (NSP) and (NS), including the Cauchy
problem of the time dependent equations we refer e.g. \cite{fei,
lio, nov, xin}. Here we prove a Liouville type  theorem of the
stationary systems (NSP)  and (NS) as follows.
 \begin{theorem} Let $N\geq 2$.
 \begin{itemize}
 \item[(i)] (Navier-Stokes Poisson equations): Suppose  $(\rho,v, \Phi)$ is a
 smooth  solution to $(NSP)$ with $k\neq0$, satisfying
 \bb\label{15}
 \|\rho\|_{L^\infty} + \|\Phi \|_{L^\infty} +\|\nabla v\|_{L^2} + \|v\|_{L^{\frac{N}{N-1}}} <\infty,
 \ee
 and the additional condition
 \bb\label{15a}
 v\in L^{\frac{3N}{N-1}} (\Bbb
 R^N)
 \ee
if $N\geq 7$. Then, $v=0$ on $\Bbb R^N$. The functions $\Phi $ and
$\rho$ are determined by solving the system:
 \bb\label{15aa} \left\{\aligned &\rho\nabla \Phi=\frac{a}{k}\nabla \rho^\gamma,\\
  &\Delta \Phi=
 \rho-\rho_0.\endaligned \right.
 \ee
\item[(ii)](Navier-Stokes equations): Suppose  $(\rho,v)$ is a
 smooth  solution to $(NS)$ satisfying
 \bb\label{15ab}
 \|\rho\|_{L^\infty} +\|\nabla v\|_{L^2} + \|v\|_{L^{\frac{N}{N-1}}} <\infty,
 \ee
 and the additional condition (\ref{15a}) if $N\geq 7$.  Then, $v=0$ and $\rho=$constant on $\Bbb
 R^N$.
\end{itemize}
\end{theorem}
{\em Remark 1.1 } Let us compare the above theorem with the previous
result in \cite{cha3}, which says that a solution to (NS)
satisfying
 \bb\label{15abc}\|\sqrt{\rho} v\|_{L^2}+ \|v\|_{L^{\frac{N}{N-1}}}
 +\|\rho^\gamma\|_{L^1}
 <\infty
 \ee
is vacuum, $\rho=0$. The main difference between (\ref{15abc}) and
(\ref{15ab}) is that in (\ref{15abc}) we are assuming the spatial
infinity is vacuum,
 \bb\label{15b}
  \rho(x)\to 0\quad\mbox{as}\quad
|x|\to \infty,
 \ee
  while in (\ref{15ab}) it is allowed to have
 \bb\label{15c}
  \rho(x)\to \rho_\infty \quad\mbox{as}\quad  |x|\to \infty
   \ee
     for a
constant $\rho_\infty>0$. Even with such non-vacuum condition at
infinity the Liouville type result holds. In particular for the 2D
stationary  Navier-Stokes equations
 if $v\in H^1(\Bbb R^2)$ and $\rho$ satisfies (\ref{15c}), then
 $v=0$ and $\rho=\rho_\infty$ on $\Bbb R^2$.\\

\section{Proof of Theorem 1.1 }
\setcounter{equation}{0}
 {\bf Proof of Theorem 1.1 }
 Let us
consider a radial cut-off function $\sigma\in C_0 ^\infty(\Bbb R^N)$
such that
 \bb\label{16}
   \sigma(|x|)=\left\{ \aligned
                  &1 \quad\mbox{if $|x|<1$}\\
                     &0 \quad\mbox{if $|x|>2$},
                      \endaligned \right.
 \ee
and $0\leq \sigma  (x)\leq 1$ for $1<|x|<2$.  Then, for each $R
>0$, we define
 \bb\label{17}
\s \left(\frac{|x|}{R}\right):=\s_R (|x|)\in C_0 ^\infty (\Bbb R^N).
 \ee
We multiply (\ref{12}) by $v\si (x)$, and integrate over $\Bbb R^N$,
and integrate by part to obtain,
 \bq\label{18}
 &&\mu \int_{\Bbb R^N} |\nabla v|^2\si \,dx +
 (\mu+\lambda) \int_{\Bbb R^N} (\text{div } v)^2\si \, dx\n\\
 &&\quad =-\mu \int_{\Bbb R^N} v \cdot(\nabla \si \cdot \nabla )v\, dx
 -(\mu+\lambda ) \int_{\Bbb R^N} \int_{\Bbb R^N}\text{ div  }v
 (v\cdot \nabla ) \si \, dx \n\\
 &&\qquad +k \int_{\Bbb R^N} \si \rho (v\cdot \nabla )\Phi \, dx-\int_{\Bbb R^N}\si (v \cdot \nabla )p \, dx
  -\int_{\Bbb R^N} \si v\cdot\text{div } (\rho v\otimes v) \, dx\n\\
 &&:=I_1+\cdots +I_5.
 \eq
Let us estimate $I_1, \cdots, I_5$ term by term.
 \bq\label{19}
 | I_1| &\leq& \frac{\mu}{R} \int_{\{R\leq |x|\leq 2R \}} |\nabla
 \s| |v||\nabla v|\, dx\n \\
 &\leq&\frac{\mu  \|\nabla
 \s\|_{L^\infty}}{R}  \left(\int_{\{R\leq |x|\leq 2R\}
 } \, dx\right)^{\frac{1}{N}}\|v\|_{L^{\frac{2N}{N-2}}(R\leq |x|\leq 2R)} \|\nabla
 v\|_{L^2(R\leq |x|\leq 2R)}\n\\
 &\leq &C \|\nabla
 \s\|_{L^\infty} \|\nabla v\|_{L^2}  \|\nabla
 v\|_{L^2(R\leq |x|\leq 2R)} \to 0
\eq
 as $R\to \infty$, where we used the Sobolev inequality,
 $\|v\|_{L^{\frac{2N}{N-2}}} \leq C \|\nabla v\|_{L^2}$.
Estimate of $I_2$ is similar to $I_1$ and we have
 \bq\label{110}
 | I_2| &\leq& \frac{\mu +\lambda }{R} \int_{\{R\leq |x|\leq 2R\} } |\nabla
 \s| |v||\nabla v|\, dx\n \\
 &\leq &C \|\nabla
 \s\|_{L^\infty} \|\nabla v\|_{L^2}  \|\nabla
 v\|_{L^2(R\leq |x|\leq 2R)} \to 0
\eq
 as $R\to \infty$.
 In order to estimate $I_3$ we  first integrate by part to obtain
  \bqn
   I_3&=& -k \int_{\Bbb R^N} \si \mathrm{div}\, (\rho v) \Phi\, dx
  -k\int_{\Bbb R^N} \Phi \rho (v\cdot \nabla)\si \, dx\n\\
  &&=-k\int_{\Bbb R^N} \Phi \rho (v\cdot \nabla)\si \, dx,
 \eqn
 where we used (\ref{11}). Therefore,
 \bq\label{111}
 |I_3|&\leq& \frac{k}{R} \|\nabla \s \|_{L^\infty}  \|\Phi\|_{L^\infty} \|\rho\|_{L^\infty}
 \int_{\{R\leq |x|\leq 2R\}} |v| \, dx \n \\
 &\leq& C\|\nabla \s \|_{L^\infty} \|\Phi\|_{L^\infty}
 \|\rho\|_{L^\infty} \|v\|_{L^{\frac{N}{N-1}}(R\leq |x|\leq 2R)}\to 0
 \eq
 as $R\to \infty$. In order to estimate $I_4$ we write the pressure
 term in the following form:
 $$
 \nabla p =a \nabla \rho^\gamma = \frac{a\gamma}{\gamma-1} \rho
 \nabla \rho^{\gamma-1}.
 $$
 Then, we have
 \bqn
  I_4&=& -\frac{a\gamma}{\gamma-1}  \int_{\Bbb R^N} \si \rho v \cdot \nabla
  \rho^{\gamma-1}\, dx\n \\
  &=&\frac{a\gamma}{\gamma-1}  \int_{\Bbb R^N} \si\,\mathrm{ div }\,
  (\rho v) \rho^{\gamma-1}\, dx +\frac{a\gamma}{\gamma-1}  \int_{\Bbb
  R^N} \rho ^{\gamma} (v \cdot \nabla ) \si \, dx\n\\
 &=&\frac{a\gamma}{\gamma-1}  \int_{\Bbb
  R^N} \rho ^{\gamma} (v \cdot \nabla ) \si \, dx
  \eqn
  thanks to (\ref{11}).
  Thus, we estimate
  \bq\label{112}
  | I_4 |&\leq & \frac{a\gamma}{(\gamma-1)R}
  \|\rho\|_{L^\infty}^\gamma \int_{\{R\leq |x|\leq 2R\}} |v||\nabla \s | \, dx
  \n \\
  &\leq & C \|\nabla \s \|_{L^\infty} \|\rho\|_{L^\infty}^\gamma
   \|v\|_{L^{\frac{N}{N-1}}(R\leq |x|\leq 2R)}\to 0
   \eq
 as $R\to \infty$. For the term $I_5$ we first compute by
 integration by part
 \bqn
 I_5&=&\int_{\Bbb R^N} \si \,|v|^2\mathrm{ div} \,(\rho v) \, dx
+\frac12 \sum_{i,j=1}^N\int_{\Bbb R^N}\si  \rho (v\cdot \nabla)
|v|^2\, dx \\
&=& -\frac12 \int_{\Bbb R^N}\si\, \mathrm{div}\, (\rho v) |v|^2   \,
dx -\frac12 \int_{\Bbb R^N}|v|^2 ( v\cdot \nabla )\si \, dx\\
&=&  -\frac12 \int_{\Bbb R^N}|v|^2 ( v\cdot \nabla )\si \, dx.
  \eqn
  Therefore, we estimate
  \bq\label{113}
  |I_5|&\leq& \frac{1}{2R}\int_{\{R\leq |x|\leq 2R\}} |v|^3 |\nabla \s|
  \,dx.
  \eq
  We estimate the right hand side of (\ref{113}), depending on $N$.
 \begin{itemize}
 \item[(i)] For $N=2$ we use the following Ladyzenskaya's inequality(\cite{lad}),
 $$
 \|v\|_{L^4}\leq 2\|v\|_{L^2}^{\frac12} \|\nabla v\|_{L^2} ^{\frac12}
 $$
to estimate
 \bq\label{114}
  |I_5|&\leq& \frac{ \|\nabla \s \|_{L^\infty} }{2R}\left( \int_{\Bbb R^2}
  |v|^4\, dx \right)^{\frac34} \left(\int_{\{ R\leq |x|\leq 2R\}} \, dx
  \right)^{\frac14}\n \\
  &\leq& C\|\nabla \s \|_{L^\infty} \|v\|_{L^4}^3 R^{-\frac12}\n \\
  &\leq& C\|\nabla \s \|_{L^\infty} \|v\|_{L^2}^{\frac32} \|\nabla v\|_{L^2} ^{\frac32} R^{-\frac12}
  \to 0
 \eq
  as $R\to \infty$.
  \item[(ii)] For $N=3$ we use the $L^p-$interpolation followed by the Sobolev inequality
 \bq\label{115}
 |I_5|&\leq&\frac{ \|\nabla \s \|_{L^\infty} }{2R}\left( \int_{\{ R\leq |x|\leq 2R\}}
  |v|^{\frac92}\, dx \right)^{\frac23}\left(\int_{\{ R\leq |x|\leq 2R\}} \, dx
  \right)^{\frac13}\n \\
  &\leq & C\|\nabla \s \|_{L^\infty}\|v\|_{L^{\frac92}(R\leq |x|\leq
  2R)}^3\n \\
  &\leq&
  C \|\nabla \s\|_{L^\infty}\|v\|_{L^{\frac32}(R\leq |x|\leq
  2R)}^{\frac13}
  \|v\|_{L^6}^{\frac83}\n \\
  &\leq&  C \|\nabla \s\|_{L^\infty}\|v\|_{L^{\frac32}(R\leq |x|\leq
  2R)}^{\frac13}
  \|\nabla v\|_{L^2}^{\frac83} \to 0
 \eq
as $R\to \infty$.
  \item[(iii)] For $4\leq N\leq 6 $ we estimate
  \bq\label{116}
  |I_5|&\leq& \frac{\|\nabla \s \|_{L^\infty}}{2R} \left(\int_{\{R\leq |x|\leq 2R\}}
  |v|^{\frac{2N}{N-2}} \, dx \right)^{\frac{3(N-2)}{2N}} \left(\int_{\{R\leq |x|\leq
  2R\}} 1 \, dx \right)^{\frac{6-N}{2N}}\n \\
  &\leq & C \|\nabla \s \|_{L^\infty}\|v\|_{L^{\frac{2N}{N-2}} (R\leq |x|\leq
  2R)} ^3 R ^{2-\frac{N}{2}}\to 0
  \eq
  as $R\to \infty$, since $\|v\|_{L^{\frac{2N}{N-2}}}\leq C\|\nabla
  v\|_{L^2}<\infty$.
 \item[(iv)] For $N\geq 7$ we use the additional hypothesis $v\in L^{\frac{3N}{N-1}} (\Bbb
 R^N)$.
  \bq\label{117}
|I_5|&\leq & \frac{\|\nabla \s \|_{L^\infty}}{2R}
\left(\int_{\{R\leq |x|\leq 2R\}} |v|^{\frac{3N}{N-1}} \, dx
\right)^{\frac{N-1}{N}}\left(\int_{\{R\leq |x|\leq
  2R\}} 1 \, dx \right)^{\frac{1}{N}} \n \\
  &\leq & C  \|\nabla \s \|_{L^\infty} \|v\|_{L^{\frac{3N}{N-1}} ( R\leq |x|\leq 2R)} ^3  \to 0
  \eq
  as $R\to \infty$.
  \end{itemize}

   Thus passing $R\to \infty$ in (\ref{18}), we find from
  (\ref{19})-(\ref{117}) that
  $$
  \mu \int_{\Bbb R^N} |\nabla v|^2 \,dx +
 (\mu+\lambda) \int_{\Bbb R^N} (\text{div } v)^2\, dx=0,
 $$
 and $v=$ constant vector in $\Bbb R^N$, which combined with the
 integrability conditions provides us with $v=0$. Thus the equation
 (\ref{12}) is reduced to the first part of (\ref{15aa}). In the
 case of (ii) the reduced equation is $\nabla p=0$, which implies
 that $\rho=$constant. $\square$\\
$$\mbox{\bf Acknowledgements}$$
 The author would like to thank Prof. Z. Xin, whose
 question for the validity of the  Liouville type result for (NS) with the non-vacuum
 boundary condition at infinity motivated this work.
 This research was supported partially by  NRF
Grant no. 2006-0093854.

\end{document}